# RESEARCH REPORT

# HOMOTOPY LIE GROUPS

JESPER M. MØLLER

ABSTRACT. Homotopy Lie groups, recently invented by W.G. Dwyer and C.W. Wilkerson [13], represent the culmination of a long evolution. The basic philosophy behind the process was formulated almost 25 years ago by Rector [32, 33] in his vision of a homotopy theoretic incarnation of Lie group theory. What was then technically impossible has now become feasible thanks to modern advances such as Miller's proof of the Sullivan conjecture [25] and Lannes's division functors [22]. Today, with Dwyer and Wilkerson's implementation of Rector's vision, the tantalizing classification theorem seems to be within grasp.

Supported by motivating examples and clarifying exercises, this guide quickly leads, without ignoring the context or the proof strategy, from classical finite loop spaces to the important definitions and striking results of this new theory.

## 1. INTRODUCTION

The aim of this report is to advertise the discovery by W.G. Dwyer and C.W. Wilkerson of a remarkable class of spaces called homotopy Lie groups or $p$-compact groups. These purely homotopy theoretic objects capture the essence of the idea of a Lie group. I shall here focus on [13], where Dwyer and Wilkerson introduce homotopy Lie groups and prove (1.1) and (1.3). Subsequent developments [11, 12, 14], possibly leading to a classification theorem in the near future, will be briefly described in the final section.

A *p-compact group* is a pointed topological space $BX$ with all of its homotopy theory concentrated at the prime $p$, whose loop space $X = \Omega BX$ satisfies a cohomological finiteness condition. The $p$-compact $r$-torus $BT = K(\mathbb{Z}_p, 2)^r$ is an example of a $p$-compact group.

A *maximal torus* for the $p$-compact group $BX$ is a map $Bi\colon BT \to BX$ of a $p$-compact torus $BT$ to $BX$ satisfying an injectivity and a maximality condition.

**Theorem 1.1** [13, 8.13, 9.4]. *Any $p$-compact group admits a maximal torus, unique up to conjugacy.*

The *Weyl group* $W_T(X)$ associated to the maximal torus $Bi\colon BT \to BX$ is, when $X$ is connected, faithfully represented in the vector space $H_2(BT;\mathbb{Q}_p)$, and the algebra homomorphism induced by the $W_T(X)$-invariant map $Bi$,

$$(1.2) \qquad H^*(Bi;\mathbb{Q}_p)\colon H^*(BX;\mathbb{Q}_p) \to H^*(BT;\mathbb{Q}_p)^{W_T(X)},$$

takes the $p$-adic rational cohomology of $BX$ into the invariant ring of $W_T(X)$.

**Theorem 1.3** [13, 9.7]. *Let $X$ be a connected $p$-compact group with maximal torus $T \to X$. Then:*

(1) *$T$ and $X$ have the same rank.*
(2) *The Weyl group $W_T(X)$ is faithfully represented as a reflection group in the $\mathbb{Q}_p$-vector space $H_2(BT;\mathbb{Q}_p)$.*
(3) *The homomorphism (1.2) is an isomorphism.*

The precise definition of a $p$-compact group, which is motivated by the concept of a finite loop space, follows shortly.

1.4. **Finite loop spaces.** Suppose $G$ is a compact Lie group. Take a free, contractible $G$-space $EG$ and define $BG = EG/G$ to be the orbit space. Then the associated fibre sequence

$$\Omega EG \to \Omega BG \to G \to EG \to BG$$

contains a homotopy equivalence $\Omega BG \to G$.

This phenomenon is embedded in the general concept of a *finite loop space*.

**Definition 1.5.** A finite loop space is a connected, pointed space $BX$ such that $X = \Omega BX$ is homotopy equivalent to a finite CW-complex.

Note that $X$ — *by definition* — is the loop space of $BX$. It is customary, though ambiguous, to refer to the finite loop space $BX$ by its underlying space $X$ and then call $BX$ the *classifying space* of $X$.

We have already seen that compact Lie groups are finite loop spaces. The classifying space of $\mathrm{SU}(2)$, for instance, is the infinite quaternionic projective space $B\,\mathrm{SU}(2) = \mathbb{H}P^\infty$. However, the class of finite loop spaces is much larger. A striking example was provided by Rector [31], who found an uncountable family of homotopically distinct finite loop spaces $BX$ with $X$ homotopy equivalent to $\mathrm{SU}(2)$. In other words, the homotopy type $\mathrm{SU}(2)$ supports uncountably many distinct loop space structures.

Rector's example destroyed all hopes of a classification theorem in the spirit of compact Lie groups — as long as one sticks to integral spaces, that is. The situation looks brighter in the category of $\mathbb{F}_p$-local spaces.

1.6. **Notation.** In the following, $p$ denotes a fixed prime number, $\mathbb{F}_p$ the field with $p$ elements, $\mathbb{Z}_p$ the ring of $p$-adic integers, and $\mathbb{Q}_p = \mathbb{Z}_p \otimes \mathbb{Q}$ the field of $p$-adic numbers.

$H^*(-)$ denotes singular cohomology with $\mathbb{F}_p$-coefficients, $H^*(-;\mathbb{F}_p)$, while $H^*(-;\mathbb{Q}_p)$ denotes $H^*(-;\mathbb{Z}_p) \otimes \mathbb{Q}$ (and *not* singular cohomology with $\mathbb{Q}_p$-coefficients).

A space $K$ is $\mathbb{F}_p$-*finite* if $H^*(K)$ is finite dimensional over $\mathbb{F}_p$. A map $A \to B$ is an $\mathbb{F}_p$-*equivalence* if it induces an isomorphism $H^*(B) \to H^*(A)$ on $H^*(-)$.

1.7. **$\mathbb{F}_p$-local spaces.** A space $K$ is $\mathbb{F}_p$-*local* if any $\mathbb{F}_p$-equivalence $A \to B$ induces a homotopy equivalence $\mathrm{map}(B, K) \to \mathrm{map}(A, K)$ of mapping spaces.

$\mathbb{F}_p$-local spaces exist — the classifying space $K(\mathbb{F}_p, i)$ for the functor $H^i(-)$ is an obvious example. In fact, any space can be made $\mathbb{F}_p$-local in a minimal and functorial way.

**Theorem 1.8** (Bousfield [4, 3.2]). *There exists a functor $K \rightsquigarrow K_p$, of the (homotopy) category of CW-complexes into itself, with a natural transformation $\eta_K \colon K \to K_p$ such that $\eta_K$ is an $\mathbb{F}_p$-equivalence and $K_p$ is $\mathbb{F}_p$-local.*

Note these categorical consequences of the definition and the theorem:

- Any $\mathbb{F}_p$-equivalence between $\mathbb{F}_p$-local spaces is a homotopy equivalence.
- A map is an $\mathbb{F}_p$-equivalence if and only if its $\mathbb{F}_p$-localization is a homotopy equivalence.
- $K$ is $\mathbb{F}_p$-local if and only if $\eta_K \colon K \to K_p$ is a homotopy equivalence.

If $K$ is nilpotent or is connected and either $H_1(K; \mathbb{F}_p) = 0$ or $\pi_1(K)$ is finite, then [5, VI.5.3, VII.3.2, VII.5.1] [13, §11] the Bousfield localization $K_p$ coincides with the (perhaps more familiar) Bousfield–Kan localization $(\mathbb{F}_p)_\infty K$. In particular [5, VI.5.2], $\pi_i(K_p) \cong \pi_i(K) \otimes \mathbb{Z}_p$ when $K$ is connected, pointed and nilpotent with finitely generated abelian homotopy groups. (This is far from true without the nilpotency hypothesis (2.2).)

1.9. **Homotopy Lie groups.** By design, actually inescapably [7], the uncountably many finite loop spaces $BX$ of Rector's example all $\mathbb{F}_p$-localize to the standard $(B\mathrm{SU}(2))_p$. This observation indicates that "$\mathbb{F}_p$-local finite loop spaces" are better behaved than integral loop spaces. The problem, however, is that this term is absurd as finite complexes are unlikely to be $\mathbb{F}_p$-local. The solution proposed in [13] is to replace the topological finiteness criterion in (1.5) by a cohomological one.

**Definition 1.10** [13, 2.2]. A *$p$-compact group* is an $\mathbb{F}_p$-local space $BX$ such that $X = \Omega BX$ is $\mathbb{F}_p$-finite.

Again, it is customary to use $X$, by definition the loop space of $BX$, when referring to the $p$-compact group $BX$, the *classifying space* of $X$.

## 2. Examples of $p$-compact groups

Let $G$ be any compact Lie group whose component group $\pi_0(G)$ is a $p$-group. Define $B\hat{G} = (BG)_p$. Then $\hat{G}$ is a $p$-compact group with $H^*(B\hat{G}) = H^*(BG)$, $\pi_0(\hat{G}) = \pi_0(G)$, and $\pi_i(\hat{G}) = \pi_i(G) \otimes \mathbb{Z}_p$ for all $i \geq 1$ [13, §11].

This example includes all finite $p$-groups such as the trivial group $\{1\}$ and the cyclic $p$-groups $\mathbb{Z}/p^n$, $n \geq 1$.

2.1. **Toral groups.** When applied to an $r$-torus $S = \mathrm{SO}(2)^r$, the above construction produces a *$p$-compact $r$-torus* $T = \hat{S}$. The classifying space $BT = K(\mathbb{Z}_p, 2)^r$ is an Eilenberg–MacLane space with homotopy in dimension two, and $H^*(BT) = \mathbb{F}_p[t_1, \ldots, t_r]$ is polynomial on $r$ generators of degree two.

Alternatively, $BT = (B\check{T})_p$, where $\check{T} = (\mathbb{Z}/p^\infty)^r$ is a *$p$-discrete $r$-torus*.

More generally, a *$p$-compact toral group* $P$ is a $p$-compact group with $BP = (B\check{P})_p$, where $\check{P}$, a *$p$-discrete toral group*, is an extension of a $p$-discrete torus $\check{T}$ by a finite $p$-group $\pi$. Note that the $\mathbb{F}_p$-localized sequence $BT \to BP \to B\pi$ is [5, II.5.1] a fibration sequence, as $\pi$ necessarily acts nilpotently on $H_i(B\check{T};\mathbb{F}_p)$.

Since $\mathrm{map}(BP,BX) \simeq \mathrm{map}(B\check{P},BX)$ for any $p$-compact group $X$ and $\check{P}$ is the union of an ascending chain of finite $p$-groups, results for finite $p$-groups can often, by *discrete approximation* [13, §6], be extended to $p$-compact toral groups.

2.2. **Exotic examples.** Call a connected $p$-compact group *exotic* if it is not of the form $\hat{G}$ for any connected compact Lie group $G$. Sullivan spheres and, more generally, (many) Clark–Ewing $p$-compact groups are exotic.

**Proposition 2.3** (Sullivan). *Assume that $n > 2$ is an integer dividing $p-1$. Then the $\mathbb{F}_p$-local sphere $(S^{2n-1})_p$ is a $p$-compact group.*

The construction goes as follows: The cyclic group $\mathbb{Z}/n$ acts on $\check{T} = \mathbb{Z}/p^\infty$ as $\mathbb{Z}/n < \mathrm{Aut}(\check{T}) \cong \mathbb{Z}_p^*$ when $n|(p-1)$. Define $BX = (B\check{N})_p$ where $\check{N} = \check{T} \rtimes \mathbb{Z}/n$ is the semidirect product. The computation
$$H^*(BX) = H^*(B\check{N}) = H^*(B\check{T})^{\mathbb{Z}/n} = \mathbb{F}_p[t]^{\mathbb{Z}/n} = \mathbb{F}_p[t^n],$$
which uses (1.7) and the fact that $n$ is prime to $p$, shows that the mod $p$ cohomology of $BX$ is polynomial on one generator of degree $2n$. Thus the $\mathbb{F}_p$-local space $BX$ is $(2n-1)$-connected [5, I.6.1], and its loop space $X$ is $(2n-2)$-connected with $H^*(X)$ abstractly isomorphic to $H^*(S^{2n-1})$. The Hurewicz theorem tells us that this abstract isomorphism is realizable by an $\mathbb{F}_p$-equivalence $S^{2n-1} \to X$, i.e. (1.7) by a homotopy equivalence $(S^{2n-1})_p \to X_p = X$.

Clark and Ewing [6] observed that applicability of Sullivan's construction is not restricted to the rank one case. Let $\check{T}$ be a $p$-discrete $r$-torus and $W < \mathrm{Aut}(\check{T}) \cong \mathrm{GL}_r(\mathbb{Z}_p)$ a finite group of order prime to $p$ acting on $\check{T}$. Put $BX = (B\check{N})_p$ where $\check{N} = \check{T} \rtimes W$. The invariant ring
$$H^*(BX) = H^*(B\check{N}) = H^*(B\check{T})^W = \mathbb{F}_p[t_1,\ldots,t_r]^W$$
is, essentially by the Shephard–Todd theorem [2, 7.2.1] [21, §23], a finitely generated polynomial algebra if and only if $W$ is a reflection group in $\mathrm{GL}_r(\mathbb{Q}_p)$. If this is the case, $H^*(X)$ is an exterior algebra on finitely many odd degree generators. In particular, $X$ is $\mathbb{F}_p$-finite and $BX$ a $p$-compact group. These Clark–Ewing $p$-compact groups are fairly well understood [8].

The Clark–Ewing $p$-compact group associated to any non-Coxeter group from the list [6] of irreducible $p$-adic reflection groups is (7.6) exotic. This scheme produces many exotics if $p$ is odd but none for $p = 2$, as the only non-Coxeter 2-adic reflection group, $W(\mathrm{DI}(4))$, number 24 on the list and of rank 3, has even order.

To come up with an exotic 2-compact group, a much more sophisticated approach is required. In their landmark paper [18], Jackowski and McClure showed how to decompose $BG$, for any compact Lie group $G$, as a generalized pushout of classifying spaces of subgroups (proper subgroups if the center of $G$ is trivial). Dwyer and Wilkerson realized that a similar decompositon applies in the case of $p$-compact groups [11, §8] and that this could be used in the construction of an exotic 2-compact group.

**Theorem 2.4** [12]. *There exists a connected $2$-compact group $\mathrm{DI}(4)$ such that $H^*(B\mathrm{DI}(4);\mathbb{F}_2)$ is isomorphic (as an algebra over the Steenrod algebra) to the* rank 4 mod 2 *Dickson algebra* $H^*(B(\mathbb{F}_2)^4 ; \mathbb{F}_2)^{\mathrm{GL}_4(\mathbb{F}_2)}$ *and* $H^*(B\,\mathrm{DI}(4) ; \mathbb{Q}_2)$ *is isomorphic to the invariant ring of* $W(\mathrm{DI}(4))$.

Assuming that such a space exists, it is possible [10] to read off from its cohomology a finite diagram that looks like the cohomological image of a diagram of spaces. Some effort is required to verify that the picture seen in cohomology actually is realizable on the level of spaces. The generalized pushout of this diagram is the exotic 2-compact group.

Since $G_1 := \mathrm{SO}(3) = \mathrm{DI}(2)$ and $G_2 = \mathrm{DI}(3)$, it is tempting to put $G_3 = \mathrm{DI}(4)$. No 2-compact group deserves the name $\mathrm{DI}(n)$ with $n \geq 5$ [19, 35].

2.5. **Cohomological invariants.** If the space $K$ is $\mathbb{F}_p$-finite, also $H^*(K ; \mathbb{Q}_p)$ is finite dimensional over $\mathbb{Q}_p$, and the *Euler characteristic*

$$\chi(K) = \sum (-1)^i \dim_{\mathbb{F}_p} H^i(K) = \sum (-1)^i \dim_{\mathbb{Q}_p} H^i(K ; \mathbb{Q}_p)$$

and the *cohomological dimension*

$$\mathrm{cd}(K) = \max\{i \mid H^i(K) \neq 0\}$$

are defined [13, 4.3, 6.13].

For a connected $p$-compact group $X$, in particular, $H^*(X ; \mathbb{Q}_p)$ is a connected finite dimensional Hopf algebra so, by Borel [3] or Milnor–Moore [26], $H^*(X ; \mathbb{Q}_p) = E(x_1, \ldots, x_r)$ is an exterior algebra on finitely many generators of odd degree, $|x_i| = 2d_i - 1$, and $H^*(BX ; \mathbb{Q}_p) = \mathbb{Q}_p[y_1, \ldots, y_r]$ is a polynomial algebra on generators of even degree, $|y_i| = 1 + |x_i|$, $1 \leq i \leq r$. The number $r = \mathrm{rk}(X)$ of generators is [13, 5.9] the *rank* of $X$. The cohomological dimension of $X$ is [14, 3.8] given by $\mathrm{cd}(X) = \max\{i \mid H^i(X ; \mathbb{Q}_p) \neq 0\} = \sum_{i=1}^r (2d_i - 1)$. For instance, $\mathrm{cd}(\hat{G}) = \dim G$, $\mathrm{rk}(\hat{G})$ is the rank of $G$, a $p$-compact $r$-torus has rank $r$, $\mathrm{cd}(P) = \mathrm{rk}(P)$ if (and only if) $P$ is a $p$-compact toral group, $\mathrm{rk}(S^{2n-1}) = 1$ and $\mathrm{cd}(S^{2n-1}) = 2n - 1$, while $\mathrm{rk}(\mathrm{DI}(4)) = 3$ and [2, Appendix A] $\mathrm{cd}(\mathrm{DI}(4)) = 27$.

**Exercise 2.6.** The trivial $p$-compact group has Euler characteristic $\chi(\{1\}) = 1$. The empty space has Euler characteristic $\chi(\emptyset) = 0$. For a connected $p$-compact group $X$, $X$ is trivial $\Leftrightarrow \chi(X) \neq 0 \Leftrightarrow \mathrm{rk}(X) = 0$ [13, 5.10].

## 3. Morphisms

A *$p$-compact group morphism* $f \colon X \to Y$ is a based map $Bf \colon BX \to BY$ between the classifying spaces. The *trivial morphism* $0 \colon X \to Y$ is the constant map $B0 \colon BX \to BY$, and the *identity morphism* $1 \colon X \to X$ is the identity map $B1 \colon BX \to BX$.

Note the fibration sequence

(3.1) $$X \xrightarrow{f} Y \to Y/f \to BX \xrightarrow{Bf} BY$$

where $Y/f$, or $Y/X$ when $f$ is understood, denotes the homotopy fibre of $Bf$.

Two morphisms $f$, $g \colon X \to Y$ are *conjugate* if the maps $Bf$, $Bg \colon BX \to BY$ are freely homotopic and $\mathrm{Rep}(X, Y) = \pi_0 \mathrm{map}(BX, BY) = [BX, BY]$ denotes the set of conjugacy classes of morphisms of $X$ to $Y$.

3.2. **Monomorphisms, epimorphisms, and isomorphisms.** The morphism $f \colon X \to Y$ is a *monomorphism* if $Y/X$ is $\mathbb{F}_p$-finite, an *epimorphism* if $Y/X$ is the classifying space of some $p$-compact group, and an *isomorphism* if $Y/X$ is contractible.

**Example 3.3.** $\{1\} \to X$ is a monomorphism with $X/\{1\} = X$, $X \to \{1\}$ is an epimorphism with $\{1\}/X = BX$, $1\colon X \to X$ is an isomorphism with $X/X = \{1\}$, and the diagonal $\Delta\colon X \to X^n$ is a monomorphism since $X^n/X$ is homotopy equivalent to $X^{n-1}$. There exists [12, 1.8] a monomorphism of $\widehat{\mathrm{Spin}(7)}$ to $\mathrm{DI}(4)$.

These definitions are motivated by

**Example 3.4.** Let $f\colon G \to H$ be a monomorphism (an epimorphism) of compact Lie groups. The homotopy fibre of the induced map $Bf\colon BG \to BH$ is $H/f(G)$ ($B(\ker f)$) so the corresponding $p$-compact group morphism $\hat f\colon \hat G \to \hat H$ is a monomorphism (an epimorphism). (Not all morphisms between $\hat G$ and $\hat H$ are induced from homomorphisms between $G$ and $H$ [16].)

A diagram $X \to Y \to Z$ of $p$-compact group morphisms is a *short exact sequence* if $BX \to BY \to BZ$ is a fibration sequence. Any $p$-compact group sits in a short exact sequence of the form $X_0 \to X \to \pi_0(X)$ where $X_0$ is the *identity component* of $X$; the identity component of a $p$-compact toral group, for instance, is a $p$-compact torus (2.1).

**Exercise 3.5.** Let $f\colon X \to Y$ and $g\colon Y \to Z$ be morphisms.
  (1) If $f$ and $g$ are monomorphisms, then $g \circ f$ is a monomorphism.
  (2) If $X$ is a $p$-compact toral group and $g \circ f$ a monomorphism, then $f$ is a monomorphism.
  (3) Assume that $X \xrightarrow{f} Y \xrightarrow{g} Z$ is a short exact sequence. Show that $f$ is a monomorphism and $g$ an epimorphism. Show also that if $X$ is a $p$-compact $r$-torus and $Z$ a $p$-compact $s$-torus, then $Y$ is a $p$-compact $(r+s)$-torus.

To be fair, part (2) of this exercise, requiring the theory of *kernels* [13, 7.1–7.3], is highly nontrivial. (The condition on $X$ can be removed [13, 9.11].)

An inspection of the Serre spectral sequence for the left segment of (3.1) yields

**Proposition 3.6** [13, 6.14]. *If $f\colon X \to Y$ is a monomorphism, then $\mathrm{cd}(Y) = \mathrm{cd}(X) + \mathrm{cd}(Y/X)$.*

**3.7. Nontrivial elements.** The existence of nontrivial elements in nontrivial $p$-compact groups represents the first and decisive step in constructing the maximal torus.

**Theorem 3.8** [13, 5.4, 5.5, 7.2, 7.3]. *Let $X$ be a nontrivial $p$-compact group.*
  (1) *There exists a monomorphism $\mathbb{Z}/p \to X$.*
  (2) *If $X$ is connected, there exists a monomorphism $S \to X$ from a $p$-compact $1$-torus $S$ to $X$.*

Note that (2) implies (1): In case the identity component $X_0$ is nontrivial, use (2) to get (3.5) a monomorphism $\mathbb{Z}/p \to S \to X_0 \to X$. Otherwise, (1) reduces to obstruction theory. (A sketch of the proof of (3.8) will be presented in §6.)

Analogously, any nontrivial, connected compact Lie group contains a copy of SO(2) (in its maximal torus).

**Exercise 3.9.** Use (3.8) and Lannes theory [22] to show that $BX$ is $\mathbb{F}_p$-finite only if $X$ is trivial. Next show that a $p$-compact group morphism which is both a monomorphism and an epimorphism is an isomorphism.

## 4. Homotopy fixed point spaces

Let $\pi$ be a finite $p$-group and $K$ a space. A $\pi$-*space* with underlying space $K$ is a fibration $K_{h\pi} \to B\pi$ over $B\pi$ with fibre $K$. A $\pi$-*map* is a map $u_{h\pi} \colon K_{h\pi} \to L_{h\pi}$ over $B\pi$.

The *homotopy orbit space* is the total space, $K_{h\pi}$, and the *homotopy fixed point space*, $K^{h\pi}$, is the space of sections (which may very well be empty). These spaces are related by the evaluation map $B\pi \times K^{h\pi} \to K_{h\pi}$.

For brevity, a $\pi$-space will often be denoted by its underlying space and a $\pi$-map by its restriction to the underlying spaces.

**Example 4.1.** The *trivial* $\pi$-space with underlying space $K$ is the trivial fibration $K \times B\pi \to B\pi$ with homotopy orbit space $K_{h\pi} = B\pi \times K$ and homotopy fixed point space $K^{h\pi} = \operatorname{map}(B\pi, K)$.

The homotopy fixed point construction $K^{h\pi}$ is functorial in both variables:

- For any $\pi$-map $u \colon K \to L$, composition with $u_{h\pi} \colon K_{h\pi} \to L_{h\pi}$ determines a map $u^{h\pi} \colon K^{h\pi} \to L^{h\pi}$.
- For any subgroup $\kappa < \pi$, any $\pi$-space is also a $\kappa$-space. The inclusion $\iota \colon \kappa \to \pi$ induces a map $K_{h\iota} \colon K_{h\kappa} \to K_{h\pi}$ over $B\iota \colon B\pi \to B\kappa$ and a map $K^{h\iota} \colon K^{h\pi} \to K^{h\kappa}$ of homotopy fixed point spaces.

The homotopy orbit space and the homotopy fixed point space are homotopy invariant constructions in that any $\pi$-map $u \colon K \to L$, which is an ordinary (nonequivariant) homotopy equivalence, induces homotopy equivalences $u_{h\pi} \colon K_{h\pi} \to L_{h\pi}$ and $u^{h\pi} \colon K^{h\pi} \to L^{h\pi}$.

4.2. **Exactness.** Let $U$ denote the (ordinary, nonequivariant) homotopy fibre of a $\pi$-map $u \colon K \to L$ (where $L$ is assumed to be connected) or, equivalently, the homotopy fibre of $u_{h\pi} \colon K_{h\pi} \to L_{h\pi}$. The pullback diagram

$$\begin{array}{ccc} U_{h\pi} & \longrightarrow & K_{h\pi} \\ \downarrow & & \downarrow u_{h\pi} \\ B\pi & \xrightarrow{l} & L_{h\pi} \end{array}$$

shows that any homotopy fixed point $l \in L^{h\pi}$ makes $U$ into a $\pi$-space such that

**Proposition 4.3** [13, 10.6]. $U \to K \to L$ *is a fibration sequence of $\pi$-maps between $\pi$-spaces and $U^{h\pi} \to K^{h\pi} \to L^{h\pi}$ is a fibration sequence of homotopy fixed point spaces (where $l \in L^{h\pi}$ serves as base point).*

4.4. **Exponential laws.** The exponential law in its simple form reads

(4.5) $$\operatorname{map}(K_{h\pi}, A) = \operatorname{map}(K, A)^{h\pi}$$

and in slightly more general form,

(4.6) $$K^{h\pi} = (K^{h\kappa})^{h(\pi/\kappa)}$$

as $K_{h\pi} = (K_{h\kappa})_{h(\pi/\kappa)}$ for any normal subgroup $\kappa \triangleleft \pi$.

The standing assumption made here that $\pi$ be finite is not essential; $\pi$ can (and will) be replaced by a $p$-discrete toral group or even a $p$-compact toral group.

## 5. Centralizers

Let $P$ be a $p$-compact toral group, $Y$ any $p$-compact group, and $g\colon P \to Y$ a morphism of $P$ to $Y$.

The *centralizer* of $g$, $C_Y(g)$, or $C_Y(P)$ when $g$ is understood, is the loop space of $BC_Y(g) = \mathrm{map}(BP, BY)_{Bg}$, the mapping space component containing $Bg$. Note the evaluation map $BC_Y(g) \times BP \to BY$. Base point evaluation $BC_Y(g) \to BY$ in particular provides the first nontrivial example of a monomorphism.

**Theorem 5.1** [13, 5.1, 5.2, 6.1]. *$C_Y(g)$ is a $p$-compact group, and $C_Y(g) \to Y$ is a monomorphism.*

The difficulty here is to show that $C_Y(g)$ and $Y/C_Y(g)$ are $\mathbb{F}_p$-finite spaces. (It is unknown if this remains true with $P$ replaced by a general $p$-compact group.)

5.2. **Central maps.** The morphism $g\colon P \to Y$ is said to be *central* if
  (1) $C_Y(g) \to Y$ is an isomorphism, or,
  (2) $g$ extends to a morphism $Y \times P \to Y$ which is the identity on $Y$.

These two conditions are equivalent as the adjoint of a morphism as in (2) is an inverse to the evaluation monomorphism in (1).

**Example 5.3** [15] [7, 2.5] [34, 9.6] [11, 12.5]. Let $C_G(f)$ be the centralizer of a homomorphism $f\colon \pi \to G$ of a finite $p$-group $\pi$ into a compact Lie group $G$ whose component group is a finite $p$-group. Then $\pi_0(C_G(f))$ is a $p$-group [17, A.4] and there is an isomorphism
$$\widehat{C_G(f)} \to C_{\hat{G}}\hat{f}$$
which is adjoint to the $\mathbb{F}_p$-localization of the map $BC_G(f) \times B\pi \to BG$ induced by the homomorphism $C_G(f) \times \pi \to G$. Thus $\hat{f}\colon \pi \to \hat{G}$ is a central morphism of $p$-compact groups if $f\colon \pi \to G$ is central as a homomorphism of Lie groups.

Here are two more examples of central morphisms.

**Theorem 5.4** [13, 5.3, 6.1]. *The constant morphism $0\colon P \to Y$ is central.*

This first example is an immediate consequence of the Sullivan conjecture as proved by H. Miller [25].

In contrast to the very deep Theorem 5.4, nothing more than elementary obstruction theory is needed for the second example of a central morphism.

**Lemma 5.5.** *The identity map $1\colon S \to S$ of a $p$-compact torus $S$ is central.*

In other words, $S$ is *abelian*. (See §8 for more on central maps and abelian $p$-compact groups.) In fact, *any* morphism into $S$ is central.

Consider a morphism $g\colon S \to Y$ defined on a $p$-compact torus. Identifying $S$ with the centralizer of the identity morphism, composition of maps provides a morphism $C_Y(g) \times S \to C_Y(g)$. The restriction to the second factor is a central factorization $g'\colon S \to C_Y(g)$ of $g$ through its own centralizer.

**Lemma 5.6** [13, 8.2, 8.3]. *Suppose that* $g\colon S \to Y$ *is a monomorphism of a p-compact torus* $S$ *to* $Y$. *Then there exists a short exact sequence*

$$S \xrightarrow{g'} C_Y(g) \to C_Y(g)/g'$$

*of p-compact groups such that* $S \xrightarrow{g'} C_Y(g) \to Y$ *is conjugate to* $g$.

Note that (5.6) asserts the existence of a classifying space $B(C_X(g)/g')$ for the homogeneous space $C_X(g)/g'$.

**Exercise 5.7.** Any monomorphism $S \to S$ is an isomorphism.

## 6. Algebraic Smith theory

Suppose that $\pi$ is a finite $p$-group and that the $p$-compact group classifying spaces $BX$ and $BY$ are $\pi$-spaces. Let $f\colon X \to Y$ be a monomorphism such that $Bf\colon BX \to BY$ is a $\pi$-map. Choose a base point $y \in (BY)^{h\pi}$ and equip $Y/X$ with the corresponding $\pi$-space structure (4.3) such that $Y/X \to BX \to BY$ is a fibration sequence of $\pi$-maps and $(Y/X)^{h\pi} \to (BX)^{h\pi} \to (BY)^{h\pi}$ a fibration sequence of homotopy fixed point spaces.

Algebraic Smith theory, based on work by J. Lannes and his collaborators and concerned with the cohomological properties, in particular the Euler characteristic (2.5), of the fibre $(Y/X)^{h\pi}$, can be summarized as follows.

**Theorem 6.1** [13, 4.5, 4.6, 5.7] [9] [24]. *Under the above assumptions the following hold*:

(1) $(Y/X)^{h\pi}$ *is* $\mathbb{F}_p$-*finite*.
(2) $\chi((Y/X)^{h\pi}) = \chi(Y/X) \mod p$.
(3) $\chi((Y/X)^{h\pi}) = \Lambda(Y/X\,;\,\pi)$ *if* $\pi$ *is cyclic*.

The *Lefschetz number* of (6.1(3)) is the alternating sum

$$\Lambda(Y/X\,;\,\pi) = \sum_{i=0}^{\infty} (-1)^i \operatorname{trace} H^i(\xi\,;\,\mathbb{Q}_p)$$

where $H^i(\xi\,;\,\mathbb{Q}_p)$ is the automorphism of $H^i(Y/X\,;\,\mathbb{Q}_p)$ induced by any generator $\xi$ of the cyclic group $\pi$.

The analogous Euler characteristic formulas were known to be true in classical Smith theory dealing with fixed point spaces for (reasonable) group actions on finite complexes.

I refrain from commenting on the proof of (6.1) but refer to [23] for more detailed information.

A particularly advantageous situation arises when the finite $p$-group $\pi$ can be replaced by a $p$-discrete torus (2.1) $\check{T}$.

**Corollary 6.2** [13, 4.7, 5.7] [9]. *Suppose that* $Bf\colon BX \to BY$ *is a* $\check{T}$-*map. Then*:

(1) $\chi((Y/X)^{hA}) = \chi(Y/X)$ *for any finite subgroup* $A < \check{T}$.
(2) $(Y/X)^{h\check{T}} \neq \emptyset$ *if* $\chi(Y/X) \neq 0$.

It is unknown if (1) also holds for infinite subgroups, such as $\check{T}$ itself.

As a demonstration of the power of Smith theory, I now sketch the proof of (3.8).

Let $X$ be a connected, nontrivial $p$-compact group so that (2.6) $\mathrm{rk}(X) > 0$. Consider the map

(6.3) $\quad \mathrm{map}(B\iota, BX)\colon \mathrm{map}(B\mathbb{Z}/p^{n+1}, BX) \to \mathrm{map}(B\mathbb{Z}/p^n, BX)$

induced by the inclusion $\iota\colon \mathbb{Z}/p^n \to \mathbb{Z}/p^{n+1}$, $n \geq 0$. Viewing $BX$ as a trivial $\mathbb{Z}/p^{n+1}$-space and $B\mathbb{Z}/p^n = (\mathbb{Z}/p)_{h\mathbb{Z}/p^{n+1}}$ as the total space of a $p$-fold covering map of $B\mathbb{Z}/p^{n+1}$, we get (4.1, 4.5)

$$\mathrm{map}(B\mathbb{Z}/p^{n+1}, BX) = (BX)^{h\mathbb{Z}/p^{n+1}}$$
$$\mathrm{map}(B\mathbb{Z}/p^n, BX) = \mathrm{map}(\mathbb{Z}/p, BX)^{h\mathbb{Z}/p^{n+1}} = (BX^p)^{h\mathbb{Z}/p^{n+1}}$$

and $\mathrm{map}(B\iota, BX) = (B\Delta)^{h\mathbb{Z}/p^{n+1}}$ reveals itself as induced by the $p$-fold diagonal (3.3) $B\Delta\colon BX \to \mathrm{map}(\mathbb{Z}/p, BX) = BX^p$, a $\mathbb{Z}/p^{n+1}$-map with fibre $X^p/X$.

By (6.1), the homotopy fibre $(X^p/X)^{h\mathbb{Z}/p^{n+1}}$ of (6.3) is $\mathbb{F}_p$-finite, and the remarkable, but straightforward, Lefschetz number calculation [13, 5.11]

$$\chi((X^p/X)^{h\mathbb{Z}/p^{n+1}}) = \Lambda(X^p/X\,;\,\mathbb{Z}/p^{n+1}) = p^{\mathrm{rk}(X)}$$

shows (2.6) that it is nonempty and noncontractible.

If all maps of $B\mathbb{Z}/p$ to $BX$ were inessential, then (5.4) $\mathrm{map}(B\iota, BX)$ would be a homotopy equivalence of the total space $\mathrm{map}(B\mathbb{Z}/p, BX) = BC_X(0)$ to the base space $BX$ of (6.3) with $n = 0$ and the homotopy fibre would be contractible — but it is not, so there must exist an essential map $Bf_1\colon B\mathbb{Z}/p \to BX$, i.e. [13, §7] a monomorphism $f_1\colon \mathbb{Z}/p \to X$. Since the homotopy fibre over $Bf_1$ of (6.3) with $n = 1$ is nonempty, $f_1$ extends (up to conjugacy) to a monomorphism [13, §7] $f_2\colon \mathbb{Z}/p^2 \to X$. Proceeding inductively, we obtain a morphism $f_\infty\colon \mathbb{Z}/p^\infty \to X$ that restricts to a monomorphism on $\mathbb{Z}/p^n$ for all $n \geq 1$. The $\mathbb{F}_p$-localization of $Bf_\infty$ is (2.1) a monomorphism $f\colon S \to X$ of a $p$-compact 1-torus $S$ to $X$.

This proves (3.8). Also (5.1) is a quick consequence of (6.1).

## 7. Maximal tori and Weyl groups

Let $X$ be any $p$-compact group. The maximal torus of $X$ is constructed by an inductive procedure.

If $X(= C_X(\{1\})/\{1\})$ is not homotopically discrete, it is (3.8) the target of a monomorphism $S_1 \to X$ defined on a $p$-compact 1-torus $S_1$. This monomorphism factors through its own centralizer (5.6) to give a short exact sequence

$$S_1 \to C_X(S_1) \to C_X(S_1)/S_1$$

of $p$-compact groups.

If $C_X(S_1)/S_1$ is not homotopically discrete, it is (3.8) the target of a monomorphism $S_2/S_1 \to C_X(S_1)/S_1$ defined on a $p$-compact 1-torus $S_2/S_1$. Pullback along this monomorphism induces a commutative diagram of $p$-compact group morphisms

$$\begin{array}{ccc} S_1 \longrightarrow S_2 \longrightarrow S_2/S_1 \\ \Big\| \quad\quad \Big\downarrow \quad\quad \Big\downarrow \\ S_1 \longrightarrow C_X(S_1) \longrightarrow C_X(S_1)/S_1 \end{array}$$

where $S_2$ is (3.5) a $p$-compact 2-torus and the middle arrow a monomorphism ($C_X(S_1)/S_2 \simeq \frac{C_X(S_1)/S_1}{S_2/S_1}$ is $\mathbb{F}_p$-finite). Thus $X$ is the target of a monomorphism (3.5) $S_2 \to C_X(S_1) \to X$ defined on a $p$-compact 2-torus.

By dimension considerations (3.6), this inductive procedure eventually stops at a *maximal torus* for $X$ where

**Definition 7.1** [13, 8.8, 8.9]. A maximal torus is a monomorphism $i\colon T \to X$ of a $p$-compact torus $T$ to $X$ such that $C_X(T)/T$ is a homotopically discrete $p$-compact group.

We have thus established the existence part of (1.1).

Let $i\colon T \to X$ be a maximal torus such that $Bi\colon BT \to BX$ is a fibration. The *Weyl space* $\mathcal{W}_T(X)$ is the topological monoid of self-maps of $BT$ over $BX$. As a space, $\mathcal{W}_T(X)$ is the fibre over $B\iota$ of $\operatorname{map}(BT, B\iota) : \operatorname{map}(BT, BT) \to \operatorname{map}(BT, BX)$, i.e. (5.5)

$$(7.2) \qquad \mathcal{W}_T(X) = (X/T)^{hT} = \coprod_w C_X(i)/C_T(w) = \coprod_w C_X(T)/T$$

where the disjoint union indexed by all $w \in \operatorname{Rep}(T, T)$, necessarily (3.5, 5.7) central automorphisms, with $i \circ w$ conjugate to $i$. The right-hand side shows that the Weyl space is homotopically discrete.

**Definition 7.3** [13, 9.6]. The *Weyl group* $W_T(X)$ is the component group (!) $\pi_0 \mathcal{W}_T(X)$ of the Weyl space.

The homotopy fixed point space $(X/T)^{hT}$ from (7.2), by discrete approximation (2.1) homotopy equivalent [13, 6.1, 6.7] to $(X/T)^{hA}$ for some finite subgroup $A < \check{T}$, is (6.1) $\mathbb{F}_p$-finite, and the computation (6.2)

$$(7.4) \qquad \chi(X/T) = \chi((X/T)^{hA}) = \chi((X/T)^{hT}) = |W_T(X)|$$

shows that the Euler characteristic of the homogeneous space $X/T$ equals the order of the Weyl group; in particular, $\chi(X/T) > 0$.

Another application of (6.2) now yields the uniqueness part of (1.1): Suppose that $i_1\colon T_1 \to X$ and $i_2\colon T_2 \to X$ are maximal tori. The fact that $(X/T_2)^{h\check{T}_1} = (X/T_2)^{hT_1} \neq \emptyset \neq (X/T_1)^{h\check{T}_2} = (X/T_1)^{hT_2}$ means that there exist morphisms $u\colon T_1 \to T_2$ and $v\colon T_2 \to T_1$, necessarily isomorphisms (3.5, 5.7), such that $i_1$ is conjugate to $i_2 \circ u$ and $i_2$ to $i_1 \circ v$.

Approaching the key results of [13], we now specialize to connected $p$-compact groups.

**Proposition 7.5** [13, 9.1]. *Let $X$ be a connected $p$-compact group and $T \to X$ a maximal torus. Then the morphism $T \to C_X(T)$ is an isomorphism.*

Again, the proof uses Smith theory.

With $C_X(T)/T = \{1\}$, (7.2) shows that the monoid morphism $W_T(X) \to \operatorname{Rep}(T, T) = [BT, BT]$ is injective, i.e. that $W_T(X)$ is faithfully represented in $H_2(BT; \mathbb{Q}_p) := H_2(BT; \mathbb{Z}_p) \otimes \mathbb{Q}$. The $W_T(X)$-invariant map $Bi\colon BT \to BX$ induces an algebra map

$$H^*(Bi; \mathbb{Q}_p) : H^*(BX; \mathbb{Q}_p) \to H^*(BT; \mathbb{Q}_p)^{W_T(X)}$$

into the invariant ring of this faithful representation. By means of a transfer map, generalizing the Becker–Gottlieb transfer [1] to fibrations with $\mathbb{F}_p$-finite fibres and ingeniously constructed [13, 9.13] via the Kan–Thurston theorem [20], Dwyer and Wilkerson show that $H^*(Bi\,;\mathbb{Q}_p)$ is injective. The extension $H^*(BX\,;\mathbb{Q}_p) \subseteq H^*(BT\,;\mathbb{Q}_p)$ is finite since the fibre $X/T$ is $\mathbb{F}_p$-finite, so these two rings have identical Krull dimensions, i.e. $\mathrm{rk}(T) = \mathrm{rk}(X)$ (1.3.(1)). The culmination of [13] is (1.3.(2)) asserting that $H^*(BX\,;\mathbb{Q}_p)$ is isomorphic to the invariant ring of the Weyl group. The fact that the invariant ring is polynomial (2.5) implies, by the classical Shephard–Todd theorem [2, 7.2.1], that (1.3.(3)) $W_T(X)$ is represented as a reflection group in the vector space $H_2(BT\,;\mathbb{Q}_p)$. Thus $W_T(X)$ must be isomorphic to a product of irreducible $p$-adic reflection groups from the Clark–Ewing list [2, 7.1], and $H^*(BX\,;\mathbb{Q}_p)$ must be isomorphic to a tensor product of the corresponding graded invariant rings.

**Example 7.6.** Suppose that $G$ is a compact Lie group with $\pi_0(G)$ a $p$-group. Then any Lie theoretic maximal torus $T \to G$ induces a maximal torus $\hat{T} \to \hat{G}$ of the $p$-compact group $\hat{G}$. The associated Weyl groups are isomorphic. The Weyl group of a $p$-compact toral group $P$ is $\pi_0(P)$. The Weyl group of the Sullivan sphere $(S^{2n-1})_p$ is $\mathbb{Z}/n$, and the Weyl group of the Clark–Ewing $p$-compact group $(B(\check{T} \rtimes W))_p$ is $W$. The Weyl group of DI(4) is $W(\mathrm{DI}(4))$, abstractly isomorphic to the product of a cyclic group of order two and the simple group of order 168.

**Exercise 7.7.** Modify the above construction of the maximal torus for a $p$-compact group to obtain an unconventional construction of the maximal torus for a compact Lie group [23] [13, 1.2].

While the $p$-adic rational cohomology $H^*(BX\,;\mathbb{Q}_p)$ is under control, it is quite another matter with the cohomology algebra $H^*(BX)$ with coefficients in $\mathbb{F}_p$. The difference between the two situations is that $H^*(BX\,;\mathbb{Q}_p)$ embeds into the polynomial ring $H^*(BT\,;\mathbb{Q}_p)$, while $H^*(BX)$ embeds into $H^*(BN_p(T))$ where $BN_p(T)$, the $p$-normalizer of the maximal torus [13, 9.8], is the Borel construction for the action of a Sylow $p$-subgroup of the Weyl group on $BT$. As a result, $H^*(BX)$ need not be polynomial but can actually be surprisingly complicated [36]. However, Dwyer and Wilkerson are able to settle perhaps the most basic structural question.

**Theorem 7.8** [13, 2.3]. *$H^*(BX)$ is a finitely generated $\mathbb{F}_p$-algebra for any $p$-compact group $X$.*

Given (7.8), it is but a small step to verify [13, 1.1] the old conjecture that $H^*(BX)$ is a finitely generated $\mathbb{F}_p$-algebra for any finite loop space $X$.

## 8. Classification

The classification scheme for $p$-compact groups, not yet completed, exhibits a pronounced analogy to Lie theory.

8.1. **Centers.** A $p$-compact group is *abelian* if its identity map is central. By (5.5), $p$-compact tori are abelian.

**Theorem 8.2** [11, 1.1] [29, 3.1]. *A $p$-compact group is abelian if and only if it is isomorphic to a product of a finite abelian $p$-group and a $p$-compact torus.*

If $Z \to X$ is a central monomorphism, then [11, 5.1] [29, 3.5] $Z$ is abelian.

**Theorem 8.3** [11, 1.2] [29, 4.4]. *For any p-compact group $X$ there exists a central monomorphism $Z(X) \to X$ such that any central monomorphism into $X$ factors, in an essentially unique way, through $Z(X)$.*

The terminal central monomorphism, essentially unique, of (8.3) is the *center* of $X$.

The (discrete approximation to the) center can be defined as the group of elements in (the discrete approximation to) $C_X(T)$ that are central in $X$. Another candidate to the center title is the centralizer of the identity morphism. Fortunately, there is no discrepancy.

**Theorem 8.4** [11, 1.3]. *The map $BZ(X) \to \mathrm{map}(BX, BX)_{B1}$, corresponding to the isomorphism $C_X(Z(X)) \to X$, is a homotopy equivalence.*

The highly nontrivial proof of (8.4) involves decomposing $BX$ as a generalized pushout — a technique also applied in the proof of (2.4).

**Example 8.5.** Let $G$ be a connected compact Lie group with Lie theoretic center $Z(G)$. Then the maps

$$(BZ(G))_p^\wedge \to \mathrm{map}(B\hat{G}, B\hat{G})_{B1} \leftarrow BZ(\hat{G})$$

adjoint to the $\mathbb{F}_p$-localization of $BZ(G) \times BG \to BG$ and $BZ(\hat{G}) \times B\hat{G} \to B\hat{G}$, respectively, are homotopy equivalences [11, 1.4, 12.1] (8.4). The center of the p-compact group $\hat{\mathrm{E}}_6$ is cyclic of order 3 if $p = 3$ and trivial otherwise.

For a connected $p$-compact group $X$ [29, 5.2],

$$(8.6) \qquad \mathrm{rk}(Z(X)) = \dim_{\mathbb{Q}_p}(\pi_1(X) \otimes \mathbb{Q}_p) = \dim_{\mathbb{Q}_p} H_2(BT; \mathbb{Q}_p)^{W_T(X)},$$

so the center is finite if and only if the fundamental group is. The quotient $p$-compact group [13, 8.3] $X/Z(X)$ has trivial center [11, 6.3] [29, 4.6].

The classification of connected $p$-compact groups essentially reduces to the simply connected case by the following result.

**Theorem 8.7** [29, 5.4]. *For any connected p-compact group $X$, there exists a short exact sequence*

$$A \to Y \times S \to X$$

*where $A$ is a finite abelian p-group, $A \to Y \times S \xrightarrow{pr_1} Y$ is a central monomorphism into the simply connected p-compact group $Y$, and $S$ is a p-compact torus.*

In (8.7), $Y$ is the universal covering $p$-compact group of $X$, $A$ is the torsion subgroup of $\pi_1(X)$, and $S = Z(X)_0$ is the identity component of the center.

**8.8. Semisimplicity.** Call a connected $p$-compact group with maximal torus $T \to X$ *simple* if the faithful representation of the Weyl group $W_T(X)$ in $H_2(BT; \mathbb{Q}_p)$ is irreducible. Sullivan spheres and Clark–Ewing $p$-compact groups (2.2) are simple by design. $\hat{G}$ is simple for any connected compact simple Lie group $G$.

For any connected $p$-compact group, the $W_T(X)$-representation $H_2(BT;\mathbb{Q}_p)$ splits as a direct sum

$$H_2(BT;\mathbb{Q}_p) = M_1 \oplus \cdots \oplus M_n$$

of irreducible representations. Provided the center $Z(X) = 0$ is trivial, all $M_i$ are nontrivial (8.6), and [14, 1.5] this splitting of $\mathbb{Q}_p[W_T(X)]$-modules descends to a splitting of $\mathbb{Z}_p[W_T(X)]$-modules,

$$H_2(BT;\mathbb{Z}_p) = L_1 \oplus \cdots \oplus L_n$$

where $L_i = H_2(BT;\mathbb{Z}_p) \cap M_i$. The splitting criterion [14, 1.4], guaranteeing the realizability as a splitting of $X$ of any splitting of the $\mathbb{Z}_p[W_T(X)]$-module $H_2(BT;\mathbb{Z}_p)$, now leads to the main result on semisimplicity.

**Theorem 8.9** [14, 1.3] [30]. *Any connected p-compact group with trivial center is isomorphic to a product of simple p-compact groups.*

For any simply connected $p$-compact group $Y$, $Z(Y)$ is finite, and the decomposition of the center free quotient $Y/Z(Y)$ into simple factors lifts to a decomposition of $Y$ into simple factors [14, 1.6].

The decompositions of (8.7) and (8.9) are also useful for the classification of endomorphisms of $p$-compact groups [28, 27].

The final stage, still unresolved, of the classification scheme has to address existence and uniqueness of simple $p$-compact groups. At the core of this problem, two obstruction groups loom.

MATEMATISK INSTITUT, UNIVERSITETSPARKEN 5, DK–2100 KØBENHAVN Ø, DENMARK
*E-mail address*: moller@math.ku.dk